\documentclass[11pt]{article}

\usepackage{amsfonts,amsmath,amscd}
 
\begin{document}
\title{Op\'erateurs diff\'erentiels globaux sur les   courbes elliptiques  }
\author{Michel Gros\thanks{membre du programme TMR de la CEE,  r\'eseau {\it{Arithmetic Algebraic
Geometry}}.}\\IRMAR, UMR CNRS 6625\\Universit\'e de Rennes I\\Campus de Beaulieu\\35042 Rennes
cedex\\France\\e-mail : michel.gros@univ-rennes1.fr}

\maketitle
\author

\setcounter{section}{-1}
\section{Introduction}

Les notes qui vont suivre sont informelles.\\

Soient $X$ une courbe projective lisse sur $k = {\Bbb{C}}$, $x \in X(k)$, $J$ la jacobienne de
$X$ et 
${\cal{D}}_{J}$ le faisceau des op\'erateurs diff\'erentiels sur celle-ci. La trivialit\'e du faisceau
tangent de $J$ implique  (cf. par exemple \cite{Rothstein}, 5.2) l' existence d' un isomorphisme
$H^{0}(J, {\cal{D}}_{J}) \simeq {\text{Sym}}_{k}\,\, H^{1}(X, O_{X})$. Soient, d' autre part,
$\pi^{0}$ l' alg\`ebre (commutative) d' Heisenberg (\cite {Fren-Ben}, example 8.4.3) (munie de sa structure
vertex quasi-conforme "canonique" sous-jacente : cf. loc. cit.) et
$H(X,x,\pi^{0} )$ l' alg\`ebre correspondante des co-invariants (cf. loc. cit. Def 8.2.7, Prop 8.4.1). On
dispose d' un isomorphisme (cf. loc. cit. example 8.4.3)) $H(X,x,\pi^{0}) \simeq {\text{Sym}}_{k}\,\,
H^{1}(X, O_{X})$ dont on d\'eduit imm\'ediatement un isomorphisme (cf. \cite {Fren-Ben}, 16.2.10 pour un
point de vue plus g\'en\'eral)\\

  $H^{0}(J, {\cal{D}}_{J}) \simeq H(X,x,\pi^{0})$.\\

Ce travail est motiv\'e par la recherche d' un  analogue de ce r\'esultat lorsque  $k$ est un
corps parfait de caract\'eristique
$p>0$ : les \'el\'ements de ${\cal{D}}_{J}$ se d\'ecrivent localement   \`a l' aide de  "puissances
divis\'ees" et il apparait naturel de substituer \`a $\pi^{0}$ la r\'eduction
${\bar{\pi}}^{0}$ modulo
$p$ de sa "forme enti\`ere" $ \pi^{0}_{ {\Bbb{Z}} }$.  \\

Nous conjecturons qu' il existe un morphisme d' alg\`ebres non trivial (qu' on esp\`ere surjectif)\\ 

$\Psi : {\bar{\pi}}^{0} \rightarrow H^{0}(J, {\cal{D}}_{J})$ \\

compatible aux morphismes de frobenius dont on dispose des deux c\^ot\'es.\\

Le probl\`eme de la d\'efinition d' un tel morphisme semble d\'ej\`a  se poser     lorsque  $X $ est une courbe
elliptique et nous donnons quelques calculs explicites allant dans ce sens.

\section{Forme   enti\`ere   de l' alg\`ebre d' Heisenberg  }

{\subsection  {L' alg\`ebre d' Heisenberg  }} 

Rappelons bri\`evement la construction de celle-ci (cf. \cite {Fren-Ben},2.1.2 et 2.1.3, 2.4 pour les
d\'etails). On dispose de l' alg\`ebre de Lie d' Heisenberg ${\cal{H}}$ qui est l' extension centrale\\

$0 \rightarrow {\Bbb{C}}.{\bf{1}} \rightarrow {\cal{H}} \rightarrow {\Bbb{C}}((t))\rightarrow 0$\\

(de cocyle $c(f,g) := - {\text{Res}}\,\,_{t=0}\,\, fdg$) et de son alg\`ebre enveloppante (compl\'et\'ee) $
{\tilde{U}}({\cal{H}})$, compl\'etion convenable de l' alg\`ebre $U({\cal{H}}')$ de g\'en\'erateurs
${\bf{1}}$ et 
$b_{n}$, $n \in {\Bbb{Z}}$ avec relations $b_{n}b_{m}-b_{m}b_{n} = n\delta_{n,-m}{\bf{1}}$,
$b_{n}{\bf{1}} - {\bf{1}}b_{n}= 0$. \\

Soient $ {\tilde{{\cal{H}}}} := {\tilde{U}}({\cal{H}})/ ({\bf{1}}-1)$ l' alg\`ebre de Weyl et 
$ {\tilde{{\cal{H}}}}_{+}$ sa sous-alg\`ebre commutative engendr\'ee par $b_{0},b_{1},...$. On peut alors
induire la repr\'esentation triviale de $ {\tilde{{\cal{H}}}}_{+}$ \`a $ {\tilde{{\cal{H}}}} $ pour obtenir
une repr\'esentation   de $ {\tilde{{\cal{H}}}} $\\

$\pi := {\text{Ind}}_{{\tilde{{\cal{H}}}}_{+}}^{{\tilde{{\cal{H}}}}}\,\, {\Bbb{C}} =
{\Bbb{C}}[b_{-1},b_{-2},...]$\\

On peut munir ce module d' une structure d' alg\`ebre vertex (que l' on peut ici pour l' essentiel ignorer) commutative quasi-conforme (cf. loc. cit. 
example 8.4.3) en faisant agir les d\'erivations
${\text{Der}}\,\, {\Bbb{C}}[[t]] = {\Bbb{C}}[[t]].\partial_{t}$  par $L_{m}:= -t^{m+1}\partial_{t}
\rightarrow \sum_{n<0} (m-n)b_{n}
\frac{\partial}{\partial b_{n-m}}$ . On note alors l' alg\`ebre (vertex) obtenue $\pi^{0}$.\\

L' alg\`ebre des coinvariants de $\pi^{0}$ (qui en est par d\'efinition un quotient)    (cf.  \cite {Fren-Ben}, 8.1.7 ) s' explicite  de la mani\`ere rappel\'ee dans l' introduction.\\

{\subsection  {Forme enti\`ere}} 

Pour la pr\'esentation de cette section, nous suivrons des id\'ees emprunt\'ees \`a \cite{Cline} (voir aussi
\cite{Shimizu}, \`a des normalisations pr\`es) qui \'eclairent sensiblement le point de vue de
\cite{Garland} pour la direction que nous souhaitons suivre. Notons d' autre part  que les vecteurs de Witt
universel qui n' apparaitront qu' en filigrane ici sont intimement li\'es (cf.
\cite{Katsura}, ¤3, \cite {Contou}, etc.) aux questions discut\'ees dans ce travailÂ.\\

On d\'efinit (cf. \cite{Cline}, ¤2 et comparer avec \cite{Mumford}, Lecture 26, B) pour $n=0,1$ des
polyn\^omes $\Lambda_{n} \in {\Bbb{Q}}[X_{1},X_{2},...]$ par l' \'egalit\'e formelle \\

$\sum_{n\geq 0} \Lambda_{n}t^{n} = {\text{exp}}\,\, (\sum_{j\geq 1} \frac{X_{j}}{j}t^{j})$.\\

{\it{Exemples}}. Ainsi, on a : \\

$\Lambda _{0} =  X_{1} $ ; $\Lambda _{1} = \frac{   [X_{2} +  X_{1}^{2} ]}{2!}$ ; $\Lambda _{2} = \frac{ [2
X_{3} + 3X_{1} X_{2} + X_{1}^{3}]}{3!}$ ; $\Lambda _{3} = \frac{  [6 X_{4} + (8X_{1}X_{3} + 3X_{2}^{2}) +
6X_{1}^{2}X_{2} + X_{1}^{4}]}{4!}$ ; $\Lambda _{4} = \frac{  [24 X_{5} + (30 X_{1}X_{4} + 20 X_{2}X_{3}) + (15
X_{1} X_{2}^{2}+ 20X_{1}^{2} X_{3}) +10 X_{1}^{3}X_{2}+ X_{1}^{5}]}{5!}$.\\

On d\'efinit maintenant un morphisme d' alg\`ebres   $\Psi_{r} : {\Bbb{Q}}[X_{1},X_{2},...]
\rightarrow
\pi^{0}$  par
$\Psi_{r}(X_{m}) := b_{rm}$ et l' on pose $\Lambda_{r}(b_{m}) := \Psi_{r}(\Lambda_{m})$.\\

{\bf{D\'efinition}} (cf. \cite {Garland}, thm. 5.8). On note ${\tilde{{\cal{H}}}}_{\Bbb{Z}}$   
  le
${\Bbb{Z}}$-r\'eseau de ${\tilde{{\cal{H}}}}$
  de base $1$ et les
$\Lambda_{r}(b_{ m})$ ($m \in {\Bbb{Z}}, r\geq 0$) et on l' appelle forme enti\`ere de
${\tilde{{\cal{H}}}}$.\\

Par le m\^eme proc\'ed\'e que pr\'ec\'edemment (ou en utilisant \cite {Garland}, 11), on d\'efinit 
 le ${\Bbb{Z}}$-r\'eseau $\pi^{0}_{\Bbb{Z}} $ de $\pi^{0}$ de base $1$ et les
$\Lambda_{r}(b_{-m})$ ($m\geq 0, r\geq 0$) et on l' appelle forme enti\`ere de $\pi^{0}$. C'est un
${\tilde{{\cal{H}}}}_{\Bbb{Z}}$ -module stable par l' action de ${\Bbb{Z}}[[t]].\partial_{t} \subset
{\text{Der}}\,\, {\Bbb{C}}[[t]]$.\\

Pour $k$ un corps parfait de caract\'eristique $p>0$, nous noterons (lorsqu' aucune confusion n' en r\'esulte)
${\bar{\pi}}^{0} :=
\pi^{0}_{\Bbb{Z}} \otimes _{\Bbb{Z}} k$.\\

\section{Op\'erateurs diff\'erentiels sur les jacobiennes de courbes}

{\subsection  {Rappels }} 

Soient $X$ un $k$-sch\'ema lisse. On dispose sur celui du faisceau ${\cal{D}}_{X}$ des op\'erateurs
diff\'erentiels au sens de EGA  . Lorsque l' on dispose de coordonn\'ees $x_{1},...,x_{n}$ sur $X$, un
\'el\'ement $P \in {\cal{D}}_{X}$ se d\'ecrit comme une somme finie $P = \sum_{{\underline{\alpha}} \in
{\Bbb{N}}^{n}} x^{\underline{\alpha}}\partial^{[{\underline{\alpha}}]}$ avec les notations et conventions
usuelles ($\partial^{[{\underline{\alpha}}]} := \partial^{[{ {\alpha}}_{1}]}...\partial^{[{ {\alpha}}_{n}]}$
; ${\alpha_{i}}! . 
\partial^{[{ {\alpha}}_{1}]}= \partial^{{ {\alpha}}_{1}}$ ; ...).\\

Le gradu\'e pour  filtration par l' ordre des op\'erateurs diff\'erentiels sur ${\cal{D}}_{X}$ est tel que\\

$gr\,\, {\cal{D}}_{X} \simeq {\text{Sym}}^{PD} {\cal{T}}_{X}$ \\

avec ${\text{Sym}}^{PD}$ l' alg\`ebre sym\'etrique "\`a puissances divis\'ees" (cf. par exemple \cite{Huyghe},
prop. 1.3.7.3 pour un \'enonc\'e impliquant ce dernier). Il ne semble toutefois pas possible en g\'en\'eral, m\^eme lorsque le fibr\'e est trivial (sch\'ema en
groupes,..) de d\'eduire de cet isomorphisme un calcul de $H^{0}(X, {\cal{D}}_{X})$.\\

{\bf{Remarque}}. La difficult\'e que l' on vient de mentionner se trouve d\'ej\`a dans le cas de $X :=
{\Bbb{G}}_{m}$ : les op\'erateurs diff\'erentiels invariants sur $X$ ne s' identifient pas \`a l' alg\`ebre
enveloppante de ${\text{Lie}}\,\, ({\Bbb{G}}_{m})$ mais \`a la "forme de Kostant" (hyperalg\`ebre) de
celle-ci.\\

On peut calculer $H^{0}(X, {\cal{D}}_{X})$ "\`a la Cech" en prenant un recouvrement de $X$ par des ouverts
affines $U_{i}$ munis de coordonn\'ees et en utilisant la description d' un \'el\'ement de $H^{0}(U_{i},
{\cal{D}}_{U_{i}})$ comme un endomorphisme de $O_{U_{i}}$ astreint aux conditions standards (\cite{EGA},
prop. 16.8.8 b)).\\

{\subsection  {La conjecture }} 

Soient $k$ un corps parfait de caract\'eristique $p>0$,  $X$ une courbe projective lisse sur $k$ et $J$ sa
jacobienne. \\

{\bf{Conjecture}}. Il existe un morphisme d' anneaux canonique non trivial (sauf si $X = {\Bbb{P}}^{1}_{k}$
!)\\

$\Psi : {\bar{\pi}}^{0} \rightarrow H^{0}(J, {\cal{D}}_{J})$. \\

Pr\'ecisons la compatibilit\'e de ce morphisme aux morphismes de frobenius dont dispose des deux
c\^ot\'es dans le cas particulier o\`u
$X$ est une courbe elliptique. On consid\`ere pour ce faire $z$ une uniformisante locale de $X$ au voisinage
de l'
\'el\'ement neutre et l' on \'ecrit, dans ce voisinage, un \'el\'ement $P \in H^{0}(J, {\cal{D}}_{J})$ sous la
forme $P =  \sum_{j\geq 0} a_{j}(z) \partial_{z}^{[j]}$ avec les $a_{j}(z)$ des s\'eries. Le morphisme de
frobenius auquel nous faisons allusion est donn\'e par $\partial_{z}^{[j]} \rightarrow
\partial_{z}^{[\frac{j}{p}]}$,
$z^{l} \rightarrow z^{\frac{l}{p}}$ (avec les conventions habituelles). Sur ${\bar{\pi}}^{0}$, le frobenius
est induit par le morphisme $\Lambda_{n} \rightarrow \Lambda_{\frac{n-1}{p}}$ (cf. \cite {Cline}, thm. 5.1)
(par exemple, si $p=2$, $  \frac{  [6 b_{-4} + (8b_{-1}b_{-3} + 3b_{-2}^{2}) +
6b_{-1}^{2}b_{-2} + b_{1}^{-4}]}{4!} \rightarrow \frac{   [b_{-2} +  b_{-1}^{2} ]}{2!}$).\\

Indiquons maintenant quels pourraient  \^etre les ingr\'edients sous-jacents \`a la construction du morphisme
$\Psi$. On dispose d' un analogue de $\Psi$ sur les vecteurs de Witt universels ${\Bbb{W}} =
{\Bbb{Z}}[[X_{1}, X_{2},...]]$\\

$\Psi^{univ} : \pi^{0}_{\Bbb{Z}} \rightarrow H^{0}({\Bbb{W}}, {\cal{D}}_{{\Bbb{W}}})$\\

induit par  
$b_{-i} \rightarrow  \frac{\partial}{\partial X_{i}} +
X_{i}\frac{\partial}{\partial X_{i+1}}+ X_{i+1}
\frac{\partial}{\partial X_{i+2}}+..$. D' autre part,    le but de $\Psi^{univ}$ est form\'e d' op\'erateurs
diff\'erentiels invariants pour la loi de groupe\footnote {$\Psi^{univ} $ est d' ailleurs
sans doute un isomorphisme de $\pi^{0}_{\Bbb{Z}} $ convenablement compl\'et\'e dans l' alg\`ebre des
op\'erateurs diff\'erentiels invariants   sur ${\Bbb{W}}$.} sur ${\Bbb{W}}$, si bien que ceux-ci
devraient se "descendre" via un morphisme se factorisant par  (cf. \cite {Katsura}, 3.10 et thm. 3.2)
$I : {\Bbb{W}} \simeq {\hat{{\Bbb{G}}}}_{a}^{\infty} = {\text{Spf}}\,\, (k[[t_{1}, t_{2},...]])
\rightarrow
{\hat{{\Bbb{G}}}}_{a}^{g} = {\text{Spf}}\,\, (k[[z_{1},...,z_{g}]])$ (avec $g$ le genre de $X$)
  d\'ecrit comme suit : on prend une base $(\omega_{0},...,\omega_{g})$ de $H^{0}(X, \Omega_{X}^{1})$,
on \'ecrit   au voisinage d' un point $P \in X$ s' envoyant sur l' \'el\'ement neutre de $J$,
$\omega_{i} = -d(\sum _{n\geq 1} I_{n}^{i} \frac{\xi^{n}}{n})$ ($\xi$ une uniformisante locale en $P$),
alors $I^{*}(z_{i}) = \sum I_{n}^{i}t_{n}$. De plus, l' invariance et la
nature de $J$ devraient suffire pour affirmer qu' on obtient ainsi des \'el\'ements de $H^{0}(J,
{\cal{D}}_{J})$.\\

\section{Exemples }

{\subsection  {Courbes elliptiques}}

C' est le premier cas non trivial. On a alors $X=J$ mais nous distinguerons parfois ces deux sch\'emas dans
les notations. Explicitons alors comment l' on construit (conjecturalement) le morphisme en question (la
compatibilit\'e aux frobenius s' obtenant par le m\^eme type de calculs).\\

On consid\`ere une \'equation de Weierstrass homog\`ene de $X$ (cf. \cite{Silver}, III, 1.)\\

$Y^{2}Z +a_{1}XYZ + a_{3}YZ^{2} = X^{3}+a_{2}X^{2}Z + a_{4}XZ^{2}+a_{6}Z^{3}$.\\

Sur un voisinage ${\cal{U}}$ de $O$ (\'el\'ement neutre pour la loi de groupe de $X$), le changement de
variables standard (cf. loc. cit. IV, 1.) conduit \`a l' \'equation\\

$w-a_{1}zw -a_{3}w^{2} = z^{3}+a_{2}z^{2}w+a_{4}zw^{2}+a_{6}w^{3}$.\\

 On introduit alors \\

$P := (a_{1}w + 3z^{2}+2a_{2}zw+a_{4}w^{2})\partial_{w} +
(1-a_{1}z-2a_{3}w-a_{2}z^{2}-2a_{4}zw-3a_{6}w^{2})\partial_{z}$\\

que l' on consid\'erera (gr\^ace, par exemple,  \`a la structure de sch\'ema en groupe de $J$) comme la
restriction \`a
${\cal{U}}$ d' un \'el\'ement de $H^{0}(J, {\cal{D}}_{J})$ not\'e encore $P$.\\

On introduit maintenant l' expression suivante de la diff\'erentielle invariante (base de $H^{0}(X,
\Omega_{X}^{1})$)\\

$\omega(z) = 1+a_{1}z+(a_{1}^{2}+a_{2})z^{2}+ (a_{1}^{3}+2a_{1}a_{2}+2a_{3})z^{3}+...)dz = (\sum_{i \geq 0}
\alpha_{i+1}z^{i}) dz$.\\

D' apr\`es ce que l' on a indiqu\'e dans le paragraphe pr\'ec\'edent, le morphisme   \\

${\bar{\pi}}^{0} \rightarrow H^{0}(J, {\cal{D}}_{J})$\\

dont je conjecture l' existence   devrait   \^etre ici induit par\\ 

$b_{-j} \rightarrow \alpha_{j}P$ ; $j\geq 1$.\\

Il est facile de v\'erifier par exemple que l' image de $\frac{b_{-1}^{2}+b_{-2}}{2!}$, c' est
\`a dire $\frac{P^{2} + a_{1}P}{2!}$ appartient bien \`a $ H^{0}(J, {\cal{D}}_{J})$. On en d\'eduit d'
ailleurs imm\'ediatement l' appartenance \`a $ H^{0}(J, {\cal{D}}_{J})$ de l' image de $\frac{b_{-2}^{2}+b_{-4}}{2!}$.\\

 Donnons maintenant quelques calculs en faveur de cette conjecture pour la famille de Legendre
$y^{2} =x(x-1)(x-\lambda)$. On a alors\\

$a_{1}=a_{3}=a_{6}=0$ ; $ a_{2} = -(1+ \lambda)$ ; $a_{4} = \lambda$.\\

$P  = ( 3z^{2}- 2(1+\lambda) zw+\lambda w^{2})\partial_{w} + (1+(1+\lambda) z^{2}-2\lambda zw)\partial_{z} =
2y \partial_{x} - (3x^{2} -2x(1+\lambda) + \lambda)\partial_{y}$\\

$\omega (z) = (1 +(-1-\lambda)z + (1+4\lambda +\lambda ^2)z^{4} + (-1-9\lambda-9\lambda^2-\lambda^3)z^{6} +
(1+16\lambda+36\lambda^{2}+16\lambda^{3}+\lambda^{4})z^{8}+...)dz$.\\

On calcule explicitement :\\

$P^{2} =   (14 \lambda x^{2} +4 \lambda^{2} x^{2}   +\lambda ^{2} +9 x^{4}-12 x^{3}+4 x^{2}-4 \lambda x  -4
\lambda^{2} x  -12 \lambda x^{3})\,\,
\partial_{y}^{2}
\, +
 \,4 y^{2} \,\,\partial_{x}^{2} \, + \, (4 \lambda y -8 \lambda xy  +12  x^{2}y -8 xy) \,\, \partial_{x}
\partial_{y} \, + \, (-4 x+6 x^{2}-4 \lambda x  +2 \lambda)\,\, \partial_{x} \, + \, (-4 y+12 xy -4 \lambda y 
)\,\, \partial_{y}$\\

$P^{3} =  (\lambda^{3}-8 x^{3}-6 \lambda^{3} x-6 x \lambda^{2}+12 \lambda x^{2}  -54 \lambda x^{5}  +12
\lambda^{3}x^{2}  -8
\lambda^{3} x^{3}  +36 \lambda^{2} x^{4}  +33 \lambda^{2} x^{2}  - 60 \lambda^{2} x^{3}  -60 \lambda x^{3}+ 99
\lambda x^{4}-54 x^{5}+36 x^{4}+27 x^{6}) \,\, \partial_{y}^{3} \,\, +\,\, 8 y^{3}\,\,\partial_{x}^{3}\,\,
+\,\, (84 \lambda x^{2}   y+24 \lambda^{2} x^{2}   y+24 y x^{2}-72\lambda y x^{3} -24
\lambda y x  +6 \lambda^{2} y + 54 x^{4} y-24 \lambda^{2} y x  -72 x^{3} y \,\, \partial_{x}
\partial_{y}^{2}\,\,  + \,\,(-24 y^{2} x-24 \lambda y^{2} x  +12 \lambda y^{2}  +36 y^{2} x^{2}) \,\, 
\partial_{x}^{2}\partial_{y}
\,\, +
 \,\, (-108 \lambda x^{2}  y+24 y x+ 84 \lambda y x - 12\lambda  y  +108 x^{3}y-12 \lambda^{2}y
+24\lambda^{2} yx -108yx^{2})\,\,
\partial_{y}^{2} \,\, +
\,\, (-72 \lambda x^{3}+54 x^{4}-24 y^{2}-72 x^{3}+6 \lambda^{2}-24 \lambda^{2} x  +24 x^{2}-24
\lambda y^{2}  +72 y^{2} x+ 84 \lambda x^{2}  +24 \lambda^{2} x^{2}  -24 \lambda x  )\,\,\partial_{x}
\partial_{y} \,\, +
 \,\,(24 y x-8 \lambda y  -8 y)\,\,\partial_{x}\, +\, (-4 \lambda-36 x^{2}+8 x+28 \lambda x -36 \lambda x^{2} 
+24 y^{2}+36 x^{3}-4 \lambda^{2}+8\lambda^{2} x)\,\,
\partial_{y} \,\, +
 \, \,(-24 \lambda y x -24 y x+12 \lambda y +36 yx^{2})\,\, \partial_{x}^{2}$\\

et l' on v\'erifie que l' on a bien  $\frac{P^{2}  }{2!}, \frac{P^{3} + 2.(-1-\lambda).P}{3!},...
\in H^{0}(J, {\cal{D}}_{J})$ ainsi que la compatibilit\'e annonc\'ee aux frobenius. \\

On peut bien s\^ur effectuer de nombreux autres calculs sur des exemples particuliers.\\

{\bf{Remarque}}; Lorsque l' on prend la courbe singuli\`ere (r\'eduction de la courbe de Tate) $y^{2}+xy =
x^{3}$, il est facile de v\'erifier que l' on a encore un morphisme \footnote {d\'ej\`a connu puisque
$X^{ns} \simeq {\Bbb{G}}_{m}$.} $\Psi : {\bar{\pi}}^{0}
\rightarrow H^{0}(X^{ns}, {\cal{D}}_{X^{ns}})$ construit sur le m\^eme principe que pr\'ec\'edemment. \\ 

{\subsection  {Op\'eration de Cartier}} 

Si ${\cal{C}}$ d\'esigne l' op\'eration de Cartier sur les formes diff\'erentielles de $X$, on sait (cf.
\cite{Stohr}, thm. 1.1, voir aussi \cite{Garnier} ) que l' on a\\

${\cal{C}}(h.\omega) = [\frac{\partial^{2p-2}}{\partial x^{p-1} \partial y^{p-1} }  (f^{p-1}h)]^{ \frac{1}{p}}.
\omega$\\

avec $f(x,y) =0$ d\'efinissant $X$ et $\omega$ la diff\'erentielle invariante.\\

Si $f(x,y) = y^{2}+a_{1}xy+a_{3}y-x^{3}-a_{2}x^{2}-a_{4}x-a_{6}$ et $p=2$, on peut v\'erifier que
$\frac{\partial^{2p-2}}{\partial x^{p-1} \partial y^{p-1} }  (f^{p-1}h) = (P+ a_{1}). h$.\\

Si $f(x,y) = y^{2} -x^{3}+(1+\lambda)x^{2} - \lambda x$ et $p=3$, on peut v\'erifier que
$\frac{\partial^{2p-2}}{\partial x^{p-1} \partial y^{p-1} }  (f^{p-1}h) = (\frac{P^{2}}{2!}+ 2.(1+\lambda)). h$.\\

Autrement dit, l' op\'erateur diff\'erentiel d\'efinissant l' op\'eration de Cartier provient bien d' un \'el\'ement
de ${\bar{\pi}}^{0}$.\\\\

\end{document}